\documentclass[12pt]{article}

\usepackage{amsmath}
\usepackage{amssymb}
\usepackage[dvips]{graphicx}

\begin{document}
\title{{\bf 
{\large Amphicheiral links with special properties, II}}
\footnotetext[0]{%
2010 {\it Mathematics Subject Classification}: 
57M25, 57M27.\\
{\it Keywords}:
amphicheiral link; Jones polynomial; HOMFLY polynomial; 
Alexander polynomial.
}}
\author{
{\footnotesize 
Teruhisa KADOKAMI}}
\date{{\footnotesize July 11, 2011}}
\maketitle

\newcommand{\circlenum}[1]{{\ooalign{%
\hfill$\scriptstyle#1$\hfill\crcr$\bigcirc$}}}

\newcommand{\svline}[1]{\multicolumn{1}{|c}{#1}}
\newfont{\bg}{cmr10 scaled\magstep4}
\newcommand{\bigzerol}{\smash{\hbox{\bg 0}}}
\newcommand{\bigzerou}{\smash{\lower1.7ex\hbox{\bg 0}}}

\newcommand{\bsquare}{\hbox{\rule{6pt}{6pt}}}
\newcommand{\qed}{\hbox{\rule[-2pt]{3pt}{6pt}}}
\newcommand{\Int}{\mathrm{Int}\ \! }
\newcommand{\Ker}{\mathrm{Ker}\ \! }
\newcommand{\Ig}{\mathrm{Im}\ \! }
\newcommand{\aug}{\mathrm{aug}\ \!}
\newcommand{\pj}{\mathrm{pr}\ \!}
\newcommand{\Tor}{\mathrm{Tor}\ \!}
\newcommand{\Spin}{\mathrm{Spin}\ \!}
\newcommand{\Eul}{\mathrm{Eul}\ \!}
\newcommand{\Vect}{\mathrm{Vect}\ \!}
\newcommand{\HULL}{\mathrm{HULL}\ \!}
\newcommand{\real}{\mathrm{real}\ \!}
\newcommand{\rank}{\mathrm{rank}\ \!}
\newcommand{\ord}{\mathrm{ord}\ \!}
\newcommand{\Sign}{\mathrm{Sign}\ \!}
\newcommand{\Hom}{\mathrm{Hom}\ \!}
\newcommand{\ad}{\mathrm{ad}\ \!}
\newcommand{\Det}{\mathrm{Det}\ \!}
\newcommand{\lk}{\mathrm{lk}\ \!}
\newcommand{\pt}{\mathrm{pt}}
\newcommand{\al}{$\alpha$}
\newcommand{\dis}{\displaystyle}

\newtheorem{df}{Definition}[section]
\newtheorem{lm}[df]{Lemma}
\newtheorem{theo}[df]{Theorem}
\newtheorem{re}[df]{Remark}
\newtheorem{pr}[df]{Proposition}
\newtheorem{ex}[df]{Example}
\newtheorem{co}[df]{Corollary}
\newtheorem{cl}[df]{Claim}
\newtheorem{qu}[df]{Question}
\newtheorem{pb}[df]{Problem}
\newtheorem{cj}[df]{Conjecture}

\makeatletter
\renewcommand{\theequation}{%
\thesection.\arabic{equation}}
\@addtoreset{equation}{section}
\makeatother

\begin{abstract}
{\footnotesize 
\setlength{\baselineskip}{10pt}
\setlength{\oddsidemargin}{0.25in}
\setlength{\evensidemargin}{0.25in}
We determine prime amphicheiral links 
with at least $2$ components and up to $11$ crossings.
There are $27$ such links.
We check also special amphicheiralities.
Most of prime links with up to $11$ crossings are detected not to be
amphicheiral by a condition on the Jones polynomial.
For the rest links, we applied conditions from the Alexander polynomial.
We added new necessary conditions for a special case.}
\end{abstract}

\section{Introduction}\label{sec:intro}
Let $L=K_1\cup \cdots \cup K_r$ be an oriented 
$r$-component link in $S^3$ with $r \ge 1$.
For an oriented knot $K$, we denote the orientation-reversed knot by $-K$.
If $\varphi$ is an orientation-reversing (orientation-preserving, respectively)
homeomorphism of $S^3$ so that
$\varphi(K_i)=\varepsilon_{\sigma(i)} K_{\sigma(i)}$ for all $i=1, \ldots, r$
where $\varepsilon_i=+$ or $-$,
and $\sigma$ is a permutation of $\{1, 2, \ldots, r\}$, 
then $L$ is said 
an {\it $(\varepsilon_1, \ldots, \varepsilon_r; \sigma)$-amphicheiral link}
(an {\it $(\varepsilon_1, \ldots, \varepsilon_r; \sigma)$-invertible link},
respectively).
A term ``amphicheiral link"  is used as a general term for
an $(\varepsilon_1, \ldots, \varepsilon_r; \sigma)$-amphicheiral link.
A link is said an {\it interchangeable link} if it is
an $(\varepsilon_1, \ldots, \varepsilon_r; \sigma)$-invertible link
such that $\sigma$ is not the identity.
An $(\varepsilon_1, \ldots, \varepsilon_r; \sigma)$-invertible link
is said an {\it invertible link} simply if there exists $1\le i\le r$
such that $\varepsilon_i=-$.
If $\sigma$ is the identity, then an amphicheiral link is said
a {\it component-preservingly amphicheiral link}, and
$\sigma$ may be omitted from the notation.
If every $\varepsilon_i=\varepsilon$ is identical for all $i=1, \ldots, r$
(including the case that $\sigma$ is not the identity), then 
an $(\varepsilon_1, \ldots, \varepsilon_r; \sigma)$-amphicheiral link
(an $(\varepsilon_1, \ldots, \varepsilon_r; \sigma)$-invertible link,
respectively)
is said an $(\varepsilon)$-amphicheiral link
(an $(\varepsilon)$-invertible link, respectively).
We use the notations $+=+1=1$ and $-=-1$.
A link $L$ with at least $2$-component is said 
an {\it algebraically split link} if the linking number of
every $2$-component sublink of $L$ is zero.
We note that 
a component-preservingly $(\varepsilon)$-amphicheiral link
is an algebraically split link.

\medskip

Necessary conditions for the Alexander polynomials of amphicheiral knots 
are studied by R.~Hartley \cite{Ha}, R.~Hartley and A.~Kawauchi \cite{HK},
and A.~Kawauchi \cite{Kw2} (cf.\ Lemma \ref{lm:HK}).
In \cite{Kw2}, non-invertibility of $8_{17}$ is firstly proved by the conditions.
On the other hand, T.~Sakai \cite{Sa} proved that
any one-variable Laurent polynomial $f(t)$ over $\mathbb{Z}$
such that $f(t)=f(t^{-1})$ and $f(1)=1$ is realized by
the Alexander polynomial of a strongly invertible knot in $S^3$.
B.~Jiang, X.~Lin, Shicheng Wang and Y.~Wu \cite{JLWW}
showed that 
(1) a twisted Whitehead doubled knot is amphicheiral
if and only if it is the unknot or the figure eight knot, and
(2) a prime link with at least $2$ components and up to $9$ crossings is 
component-preservingly $(+)$-amphicheiral if and only if it is the Borromean rings.
They used S.~Kojima and M.~Yamasaki's $\eta$-function \cite{KY}.
Shida Wang \cite{Wa} determined 
prime component-preservingly $(+)$-amphicheiral links 
with at least $2$ components and with up to $11$ crossings by the same method 
as \cite{JLWW}. 
There are four such links
(cf.\ Theorem \ref{th:table} (3)).
For geometric studies of symmetries of arborescent knots, see
F.~Bonahon and L.~C.~Siebenmann \cite{BS}.
The author \cite{Kd2} studied necessary conditions for 
the Alexander polynomials of algebraically split 
component-preservingly amphicheiral links 
by computing the Reidemeister torsions of surgered manifolds along the link
(cf.\ Lemma \ref{lm:2-comp}).
The author and A.~Kawauchi \cite{KK} obtained necessary conditions
by invariants deduced from the quadratic forms of a link \cite{Kw1, Kw5},
and by using the conditions
they showed that the Alexander polynomial of
an algebraically split component-preservingly $(\varepsilon)$-amphicheiral link
with even components is zero 
(cf.\ Conjecture \ref{cj:Kd2}) and
determined amphicheiral links with up to $9$ crossings
(cf.\ Lemma \ref{lm:sign}).

\medskip

We determine prime amphicheiral links with at least $2$ components
and up to $11$ crossings.
For a link with the crossing number up to $9$, 
we use the notation of D.~Rolfsen's book \cite{Ro}, and
for a link with the crossing number $10$ or $11$, 
we use a slightly modified notation from M.~Thistlethwaite's table on
D.~Bar-Natan and S.~Morrison's website \cite{BM}
(see Section \ref{sec:table}).
In the present paper,
we used information of the Jones polynomial and the multi-variable
Alexander polynomial for the class on the website \cite{BM}.
For prime links with up to $10$ and $11$ crossings,
firstly we checked a condition on the Jones polynomial
(cf.\ Lemma \ref{lm:Jones}).
Most of them are ruled out by the condition.
For the rest links, we applied the same conditions 
on the Alexander polynomial as in \cite{KK}
(cf.\ Lemma \ref{lm:amp}, Lemma \ref{lm:sub}, 
Lemma \ref{lm:even}, Lemma \ref{lm:2-comp}
and Lemma \ref{lm:sign}).
$10_{a51}^2$ and $11_{n127}^2$ could not be detected not to be amphicheiral
by the conditions and the HOMFLY polynomial
(see Figure 4).
Thus we made the following condition for a special case
(cf.\ Lemma \ref{lm:branch})
which is inspired by the method in \cite{JLWW}, 
and the result of \cite{HK, Kw2} on the Alexander polynomials 
of amphicheiral knots.
The conditions would be useful in determining
the {\it link-symmetric group} of a link
(cf.\ \cite{Hi, Kd1, Wh}).

\medskip

For an $r$-component link $L=K_1\cup \ldots \cup K_r$,
let ${\mit \Delta}_L(t_1, \ldots, t_r)$ be 
the $r$-variable Alexander polynomial of $L$
which is an element of the $r$-variable Laurent polynomial ring 
$\mathbb{Z}[t_1^{\pm 1}, \ldots, t_r^{\pm 1}]$
over $\mathbb{Z}$ where $t_i$\ ($i=1, \ldots, r$)
is a variable corresponding to a meridian of $K_i$.
For two elements $A$ and $B$ in 
$\mathbb{Z}[t_1^{\pm 1}, \ldots, t_r^{\pm 1}]$
($(\mathbb{Z}/2\mathbb{Z})[t_1^{\pm 1}, \ldots, t_r^{\pm 1}]$, 
respectively), we denote by 
$A\doteq B$ ($A\doteq_2 B$, respectively)
if they are equal up to multiplications of trivial units.
A one variable Laurent polynomial $r(t)\in \mathbb{Z}[t^{\pm 1}]$
is {\it of type $X$} if there are integers
$n\ge 0$ and $\lambda \ge 3$ with $\lambda$ odd, and
$f_i(t)\in \mathbb{Z}[t, t^{-1}]$
$(i=0, 1, \ldots, n)$
such that
$f_i(t)\doteq f_i(t^{-1})$, 
$|f_i(1)|=1$, and for $i>0$,
$f_i(t)\doteq_2 f_0(t)^{2^i}p_{\lambda}(t)^{2^{i-1}}$
where $p_{\lambda}(t)=(t^{\lambda}-1)/(t-1)$, and
$$r(t)\doteq 
\left\{
\begin{array}{ll}
f_0(t)^2 & (n=0),
\medskip\\
f_0(t)^2f_1(t)\cdots f_n(t) & (n\ge 1).
\end{array}
\right.$$

\begin{lm}\label{lm:branch}
Let $L=K_1\cup K_2$ be an oriented $2$-component link in $S^3$
such that the linking number $\ell$ of $L$ is not zero, and
$K_2$ is the trivial knot.
Let $\tilde{K}$ be the lifted knot of $K_1$
in the $p$-fold branched covering
over $K_2$ where $p\ge 2$ is coprime to $\ell$.
\begin{enumerate}
\item[(1)]
If $L$ is component-preservingly amphicheiral,
then $\tilde{K}$ is an amphicheiral knot.

\item[(2)]
If $L$ is component-preservingly $(-, +)$-amphicheiral,
then $\tilde{K}$ is a $(-)$-amphicheiral knot, and
there is an element $f(t)\in \mathbb{Z}[t, t^{-1}]$
such that
$|f(1)|=1$, $f(t^{-1})\doteq f(-t)$, and
$${\mit \Delta}_{\tilde{K}}(t^2)\doteq
\prod_{i=1}^{p-1}
{\mit \Delta}_L(t^2, \zeta_p^i)
\doteq f(t)f(t^{-1})$$
where $\zeta_p$ is a primitive $p$-th root of unity.

\item[(3)]
If $L$ is component-preservingly $(+, -)$-amphicheiral,
then $\tilde{K}$ is a $(+)$-amphicheiral knot, and
there are 
$r_j(t)\in \mathbb{Z}[t, t^{-1}]$ of type $X$
and an odd number $\alpha_j$\ $(j=1, \ldots, m)$
such that
$${\mit \Delta}_{\tilde{K}}(t)\doteq
\prod_{i=1}^{p-1}
{\mit \Delta}_L(t, \zeta_p^i)
\doteq \prod_{j=1}^m r_j(t^{\alpha_j})$$
where $\zeta_p$ is a primitive $p$-th root of unity.
In particular, if $\tilde{K}$ is hyperbolic, then
we can take $m=1$ and $\alpha_1=1$.

\end{enumerate}
\end{lm}

\begin{co}\label{co:branch}
Under the same setting as Lemma \ref{lm:branch},
if $L$ is $(-, +)$-amphicheiral, then
in the prime factorization of $|{\mit \Delta}_{\tilde{K}}(-1)|$,
the power of a prime which is congruent to $3$ modulo $4$ is even.
In particular, if $\ell$ is odd, then
$|{\mit \Delta}_L(-1, -1)|$ satisfies the condition.
\end{co}

Let $\mathcal{A}_n$ ($\mathcal{C}_n$, respectively)
be the set of prime amphicheiral links 
(component-preservingly amphicheiral links, respectively)
with at least $2$ components and up to $n$ crossings, and 
$\mathcal{A}_n^{\varepsilon}$ 
the subset of $\mathcal{A}_n$ consisting of
$(\varepsilon)$-amphicheiral links
($\mathcal{C}_n^{\varepsilon}$
the subset of $\mathcal{C}_n$ consisting of
component-preservingly $(\varepsilon)$-amphicheiral links, 
respectively)
where $\varepsilon=+$ or $-$.
It is clear that
$\mathcal{A}_n\supset \mathcal{C}_n$,
$\mathcal{A}_n^{\pm}\supset \mathcal{C}_n^{\pm}$,
$\mathcal{A}_n\supset \mathcal{A}_n^{\pm}$
and $\mathcal{C}_n\supset \mathcal{C}_n^{\pm}$.

\begin{theo}\label{th:table}
Under the setting above, we have the following:
\begin{enumerate}
\item[(1)]
$\mathcal{C}_{11}=\{2_1^2, 6_2^2, 6_2^3, 8_8^2, 8_6^3, 8_3^4, $
$$\begin{array}{l}
10_{a56}^2, 10_{a81}^2, 10_{a83}^2, 10_{a86}^2, 10_{a116}^2, 
10_{a120}^2, 10_{a121}^2, 10_{a136}^3, 10_{a140}^3, 10_{a169}^4, 
\medskip\\
10_{n36}^2, 10_{n46}^2, 10_{n107}^4\}.
\end{array}$$

\item[(2)]
$\mathcal{A}_{11}\setminus \mathcal{C}_{11}=
\{8_4^3, 9_{61}^2, 
10_{a151}^3, 10_{a156}^3, 10_{a158}^3, 10_{n59}^2, 10_{n105}^4, 
11_{n247}^2\}$.

\item[(3)]
$\mathcal{C}_{11}^+=\{ 6_2^3, 
10_{a140}^3, 10_{n36}^2, 10_{n107}^4\}$, 
$\mathcal{C}_{11}^-=\{10_{n36}^2\}$.

\item[(4)]
$\mathcal{A}_{11}^+\setminus \mathcal{C}_{11}^+=
\{8_4^3, 8_6^3, 8_3^4, 
10_{a151}^3, 10_{n59}^2, 11_{n247}^2\}$, 

$\mathcal{A}_{11}^-\setminus \mathcal{C}_{11}^-=
\{6_2^3, 8_4^3, 8_6^3, 8_3^4, 
10_{a140}^3, 10_{a151}^3, 10_{n59}^2, 10_{n107}^4, 
11_{n247}^2\}$.

\end{enumerate}
\end{theo}

We remark that Theorem \ref{th:table} (3)
corresponds to the theorem of Shida Wang \cite{Wa}.

\medskip

The following would also be useful if we determine
prime amphicheiral links with the crossing number greater than $11$
in the future
(see also Conjecture \ref{cj:Tait} in Section \ref{sec:remark}).

\begin{lm}\label{lm:Tait}
The minimal crossing number of an alternating amphicheiral link is even.
\end{lm}

In Section \ref{sec:proof}, we prove Lemma \ref{lm:branch}, 
Corollary \ref{co:branch} and Lemma \ref{lm:Tait}.
In Section \ref{sec:table}, we prove Theorem \ref{th:table}.
In Section \ref{sec:remark}, 
we give some remarks related to our previous results
\cite{Kd2}.

\section{Proof of Lemma \ref{lm:branch}, 
Corollary \ref{co:branch}
and Lemma \ref{lm:Tait}}\label{sec:proof}

To show Lemma \ref{lm:branch},
we need two results.

\medskip

K.~Murasugi \cite{Mu1} gave a formula of 
the Alexander polynomial of a periodic knot in $S^3$,
M.~Sakuma \cite{Sk} extended it to the case of periodic links,
and V.~G.~Turaev \cite{Tr}
extended them for more general settings.
In the present paper,
we use only the case of periodic knots in $S^3$.

\begin{lm}\label{lm:period}
{\rm (Murasugi \cite{Mu1}; Sakuma \cite{Sk}; Turaev \cite{Tr})}
Let $L=K_1\cup \ldots \cup K_r \cup K_{r+1}$ 
be an $(r+1)$-component link in $S^3$ such that
$r\ge 1$ and $K_{r+1}$ is the trivial knot.
Let $\tilde{L}=\tilde{L}_1\cup \ldots \cup \tilde{L}_r$
be the lifted link in the $p$-fold branched covering
over $K_{r+1}$
where $p\ge 2$ and $\tilde{L}_i$ is the lifted link of $K_i$
$(i=1, \ldots, r)$.
Then $\tilde{L}$ is a $p$-periodic link in $S^3$, and
we have
$${\mit \Delta}_{\tilde{L}}(t_1, \ldots, t_r)
\doteq
\prod_{i=1}^{p-1}{\mit \Delta}_L(t_1, \ldots, t_r, \zeta_p^i)$$
where $\zeta_p$ is a primitive $p$-th root of unity.
\end{lm}

R.~Hartley \cite{Ha},
R.~Hartley and A.~Kawauchi \cite{HK}, and A.~Kawauchi \cite{Kw2}
gave necessary conditions on the Alexander polynomials
of amphicheiral knots.

\begin{lm}\label{lm:HK}
{\rm (Hartley \cite{Ha}; Hartley and Kawauchi \cite{HK}; Kawauchi \cite{Kw2})}
\begin{enumerate}
\item[(1)]
Let $K$ be a $(-)$-amphicheiral knot.
Then there is an element $f(t)\in \mathbb{Z}[t, t^{-1}]$
such that
$|f(1)|=1$, $f(t^{-1})\doteq f(-t)$, and
$${\mit \Delta}_K(t^2)\doteq f(t)f(t^{-1}).$$

\item[(2)]
Let $K$ be a $(+)$-amphicheiral knot.
Then there are 
$r_j(t)\in \mathbb{Z}[t, t^{-1}]$ of type $X$
(cf.\ Section \ref{sec:intro})
and an odd number $\alpha_j$\ $(j=1, \ldots, m)$
such that
$${\mit \Delta}_K(t)\doteq \prod_{j=1}^m r_j(t^{\alpha_j}).$$
In particular, if $K$ is hyperbolic, then
we can take $m=1$ and $\alpha_1=1$.

\end{enumerate}
\end{lm}

\noindent
{\bf Proof of Lemma \ref{lm:branch}}\ 
Let $L=K_1\cup K_2$ be 
an oriented $2$-component component-preservingly amphicheiral link 
in $S^3$ with the linking number $\ell \ne 0$.
Since the exterior of $L$ is orienation-preserving homeomorphic
to that of the mirror image $L^*$ of $L$
with preserving boundary components
(which may not preserve orientations),
$\tilde{K}$ is amphicheiral for every $p$
which shows (1).
Suppose that $L$ is $(-, +)$-amphicheiral.
Then it is easy to see that
$\tilde{K}$ is a $(-)$-amphicheiral knot, and
the condition of the Alexander polynomial
is obtained from Lemma \ref{lm:period} and Lemma \ref{lm:HK} (1).
Suppose that $L$ is $(+, -)$-amphicheiral.
Then it is easy to see that
$\tilde{K}$ is a $(+)$-amphicheiral knot, and
the condition of the Alexander polynomial
is obtained from Lemma \ref{lm:period} and Lemma \ref{lm:HK} (2).
\qed

\bigskip

\noindent
{\bf Proof of Corollary \ref{co:branch}}\ 
By substituting $t=\sqrt{-1}$ in the equation
in Lemma \ref{lm:branch} (2)
(or Lemma \ref{lm:HK} (1)), 
and a result from elementary number theory
on primes in $\mathbb{Z}[\sqrt{-1}]$,
we have the result.
\qed

\bigskip

To show Lemma \ref{lm:Tait},
we need two results.

\medskip

Let $D=D_1\cup \ldots \cup D_r$ be an oriented link diagram
for an $r$-component link $L=K_1\cup \ldots \cup K_r$,
$w(D)$ the {\it writhe} of $D$ which is the sum of
the signs of the crossings, and
$c(D)$ the the {\it crossing number} of $D$.

\begin{lm}\label{lm:writhe}
In the situation above,
let $D'$ be an oriented diagram obtained from $D$
by reversing the orientation of the $i$-th component $D_i$.
Then we have
$$w(D')=w(D)-4\sum_{
{\scriptstyle 1\le j\le r}
\atop {\scriptstyle j\ne i}}\mathrm{lk}\ \! (K_i, K_j)$$
where $\mathrm{lk}\ \! (K_i, K_j)$ is the linking number
of $K_i$ and $K_j$, and
$w(D)\ (\mathrm{mod}\ \! 4)$ does not depend 
on the orientation of $D$.
\end{lm}

W.~Menasco and M.~Thistlethwaite \cite{MT}
gave the affirmative answer for Tait's flyping conjecture.
For an alternating link $L$,
a {\it reduced diagram} of $L$ is a diagram of $L$
which is an alternating diagram without nugatory crossings.
A {\it flyping} is an operation on a link diagram.

\begin{lm}\label{lm:flype}
{\rm (Menasco and Thistlethwaite \cite{MT})}
Let $L$ be an oriented prime alternating link.
Let $D$ and $D'$ be two reduced diagrams of $L$.
Then $D$ and $D'$ are related by a finite sequence
of flypings.
As consequences, we have $w(D)=w(D')$, and 
$c(D)=c(D')$ which is the minimal crossing number of $L$.
\end{lm}

\noindent
{\bf Proof of Lemma \ref{lm:Tait}}\ 
Let $L$ be an oriented prime alternating link,
and $D$ a reduced diagram of $L$.
By Lemma \ref{lm:flype},
$w(D)$ is an invariant of $L$.
By Lemma \ref{lm:writhe}, 
$w(D)\ (\mathrm{mod}\ \! 4)$
is an invariant of $L$ as an unoriented link.
Let $D^*$ be the mirror image diagram of $D$.
Suppose that $L$ is amphicheiral.
Then we have $w(D^*)=-w(D)\equiv w(D)\ (\mathrm{mod}\ \! 4)$,
and hence $2w(D)\equiv 0\ (\mathrm{mod}\ \! 4)$.
It implies that the crossing number of $L$ is even.
Since a reduced diagram for a non-prime alternating link is realized
by connected sums of reduced diagrams of the prime factors
(see \cite{Mu2}), we have the result.
\qed

\section{Amphicheiral links with up to $11$ crossings}\label{sec:table}
In this section, we determine
prime amphicheiral links with at least $2$ components and up to $11$ crossings.
For a link with the crossing number up to $9$, 
we use the notation of D.~Rolfsen's book \cite{Ro}, and
for a link with the crossing number $10$ or $11$, 
we use a slightly modified notation from M.~Thistlethwaite's table on
D.~Bar-Natan and S.~Morrison's website \cite{BM} below.
In Rolfsen's table \cite{Ro},
an $r$-component link such that $r \ge 2$
and the crossing number $c$
is denoted by $c_k^r$
where $k$ is the ordering of the link in the table.
In Thistlethwaite's table \cite{BM},
an $r$-component link such that $r \ge 2$
and the crossing number $c$
is denoted by $\mathrm{L}cak$ or $\mathrm{L}cnk$
where 
`$a$' implies that the link is alternating,
`$n$' implies that the link is non-alternating,
and $k$ is the ordering of the link in the table.
We modify the notations $\mathrm{L}cak$ and $\mathrm{L}cnk$
into $c_{ak}^r$ and $c_{nk}^r$, respectively.

\medskip

We raise some conditions on the Jones polynomial and the Alexander polynomial
of an amphicheiral link without proofs.
For an oriented $r$-component link 
$L=K_1\cup \ldots \cup K_r$,
let $V_L(t)$ be the Jones polynomial of $L$
with one variable $t$, and
$P_L(m, l)$ the HOMFLY polynomial of $L$
with two variables $m$ and $l$.
Then $V_L(t)$ is an element of $\mathbb{Z}[t^{\pm \frac 12}]$,
and $P_L(m, l)$ is an element of $\mathbb{Z}[m^{\pm 1}, l^{\pm 1}]$.
Let $L_{\varepsilon_1, \ldots, \varepsilon_r}
=\varepsilon_1 K_1\cup \ldots \cup \varepsilon_r K_r$
be an oriented link obtained from $L$
by changing the oriented $i$-th component $K_i$
into $\varepsilon_i K_i$\ ($i=1, \ldots, r$)
where $\varepsilon_i=+$ or $-$,
and $-L=L_{-, \ldots, -}$.
Let $L^*=K_1^*\cup \ldots \cup K_r^*$ 
be the mirror image of $L$ with the induced orientation.

\begin{lm}\label{lm:Jones}
Under the settings above, we have the following:
\begin{enumerate}
\item[(1)]
$V_L(t)\in t^{\frac{r+1}{2}}\cdot
\mathbb{Z}[t^{\pm 1}]$,
$V_{-L}(t)=V_L(t)$,
and $V_{L^*}(t)=V_L(t^{-1})$.

\item[(2)]
$V_{L_{\varepsilon_1, \ldots, \varepsilon_r}}(t)
=t^a\cdot V_L(t)$
where
${\displaystyle
a=\frac 32\sum_{1\le i<j\le r}
(1-\varepsilon_i\varepsilon_j)\mathrm{lk}\ \! (K_i, K_j)}$
if $r\ge 2$, and $a=0$ if $r=1$.

\item[(3)]
If $L$ is $(\varepsilon_1, \ldots, \varepsilon_r; \sigma)$-amphicheiral,
then we have
$V_L(t^{-1})=t^a\cdot V_L(t)$
where $a$ is the same as in (2).

\item[(4)]
If $L$ is amphicheiral, then
$V_L(t^{-1})$ is equal to $V_L(t)$
up to multiplication of $t^k$\ $(k\in \mathbb{Z})$
(i.e.\ the coefficients of $V_L(t)$ are symmetric).

\item[(5)]
$P_L(m, l)\in (ml)^{r+1}\cdot
\mathbb{Z}[m^{\pm 2}, l^{\pm 2}]$,
$P_{-L}(m, l)=P_L(m, l)$,
and $P_{L^*}(m, l)=P_L(m, l^{-1})$.

\item[(6)]
If $L$ is $(\varepsilon_1, \ldots, \varepsilon_r; \sigma)$-amphicheiral,
then we have
$P_{L_{\varepsilon_1, \ldots, \varepsilon_r}}(m, l)
=P_L(m, l^{-1}).$

\item[(7)]
$V_L(t)=P_L\left(\sqrt{-1}
\left(t^{\frac 12}-t^{-\frac 12}\right), \sqrt{-1}t\right)$.

\end{enumerate}
\end{lm}

Lemma \ref{lm:Jones} (1), (2), (5) and (7) are basic properties.
It is easy to see that (3) and (4) are deduced by (1) and (2), 
and (6) is deduced by (5).
Let $L_+$, $L_-$ and $L_0$ be three oriented links such that
they are identical except the local parts as in Figure 1.

\begin{figure}[htbp]
\begin{center}
\includegraphics[scale=0.7]{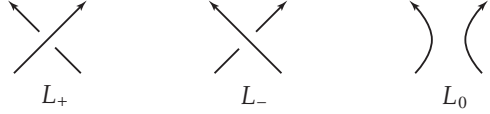} 
\label{fig:skein}
\caption{$L_+$, $L_-$ and $L_0$}
\end{center}
\end{figure}

We computed the HOMFLY polynomial of a link by the following
skein relation.
$$
lP_{L_+}(m, l)+l^{-1}P_{L_-}(m, l)+mP_{L_0}(m, l)=0,\quad
P_U(m, l)=1
$$
where $U$ implies the trivial knot.

\begin{lm}\label{lm:amp}
{\rm (\cite[Lemma 2.5]{Kd2})}
Let $L=K_1\cup \cdots \cup K_r$ be
an $r$-component 
$(\varepsilon_1, \ldots, \varepsilon_r; \sigma)$-amphicheiral link
where $\varepsilon_i=+$ or $-$\ $(i=1, \ldots, r)$, and
$\sigma$ is a permutation of $\{1, 2, \ldots, r\}$.
Then we have
$${\mit \Delta}_L(t_1, \ldots, t_r)\doteq 
{\mit \Delta}_L\left(t_{\sigma(1)}^{\varepsilon_{\sigma(1)}}, \ldots, 
t_{\sigma(r)}^{\varepsilon_{\sigma(r)}}\right).$$
\end{lm}

\begin{lm}\label{lm:sub}
{\rm (\cite[Lemma 3.1]{Kd2}, \cite[Lemma 4.2]{KK})}
Let $L=K_1\cup \cdots \cup K_r$ be 
an oriented $r$-component link.
\begin{enumerate}
\item[(1)]
If $L$ is an 
$(\varepsilon_1, \ldots, \varepsilon_r; \sigma)$-amphicheiral link, 
then a sublink $L'=K_{i_1}\cup \cdots \cup K_{i_s}$\ 
$(1\le i_1<\cdots<i_s\le r)$ is an
$(\varepsilon_{i_1}, \ldots, \varepsilon_{i_s}; \rho)$-amphicheiral link
where $\sigma$ is a permutation of $\{1, 2, \ldots, r\}$, and
$\rho$ is a permutation of $\{i_1, i_2, \ldots, i_s\}$
induced by $\sigma$.

\item[(2)]
If $r \ge 3$ is odd, and
$\ell_{1,2}\cdot \ell_{2,3}\cdots 
\ell_{r-1,r}\cdot \ell_{r,1}\ne 0$
where $\ell_{p,q}$ is the linking number of $K_p$ and $K_q$,
then $L$ is not component-preservingly amphicheiral.

\item[(3)]
If $\ell_{1,2}\cdot \ell_{2,3}\cdot \ell_{3,1}\ne 0$
where $\ell_{p,q}$ is the linking number of $K_p$ and $K_q$,
then $L$ is not amphicheiral.

\end{enumerate}
\end{lm}

\begin{lm}\label{lm:even}
Let $L=K_1\cup K_2$ be a $2$-component link
with non-zero even linking number $e$.
Then we have the following:
\begin{enumerate}
\item[(1)]
{\rm (Hartley \cite{Ha}, \cite[Lemma 3.2]{Kd2})}
$L$ is not component-preservingly amphicheiral.

\item[(2)]
{\rm (\cite[Lemma 4.3]{KK})}
If $e\equiv 2\ (\mathrm{mod}\ \! 4)$, then
$L$ is not $(\pm, \mp; (1\ 2))$-amphicheiral
where $(1\ 2)$ is the non-trivial permutation
of $\{1, 2\}$.

\end{enumerate}
\end{lm}

\begin{lm}\label{lm:2-comp}
{\rm (\cite[Corollary 1.4]{Kd2})}
If $L=K_1\cup K_2$ is
an algebraically split component-preservingly amphicheiral link, 
then ${\mit \Delta}_L(t_1, t_2)$ is divisible by
$(t_1-1)^2(t_2-1)^2$.
\end{lm}

\begin{lm}\label{lm:sign}
{\rm (\cite[Corollary 1.2]{KK})}
Let $L=K_1\cup \ldots \cup K_r$ be
an $r$-component amphicheiral link
such that $r+\ell(L)$ is even
where $\ell(L)$ is the total linking number.
Then we have
$${\mit \Delta}_L(-1, \ldots, -1)=0.$$
In particular, if $L$ is an $(\varepsilon)$-amphicheiral link
where $\varepsilon=+$ or $-$, and $r$ is even,
then we have
$${\mit \Delta}_L(t, \ldots, t)=0.$$
\end{lm}

Prime links with up to $9$ crossings
have been determined in \cite{KK} (cf.\ Figure 2).
From now on, we restrict the case that a link is prime
with the crossing number $10$ or $11$.

\newpage

\begin{figure}[htbp]
\begin{center}
\includegraphics[scale=0.6]{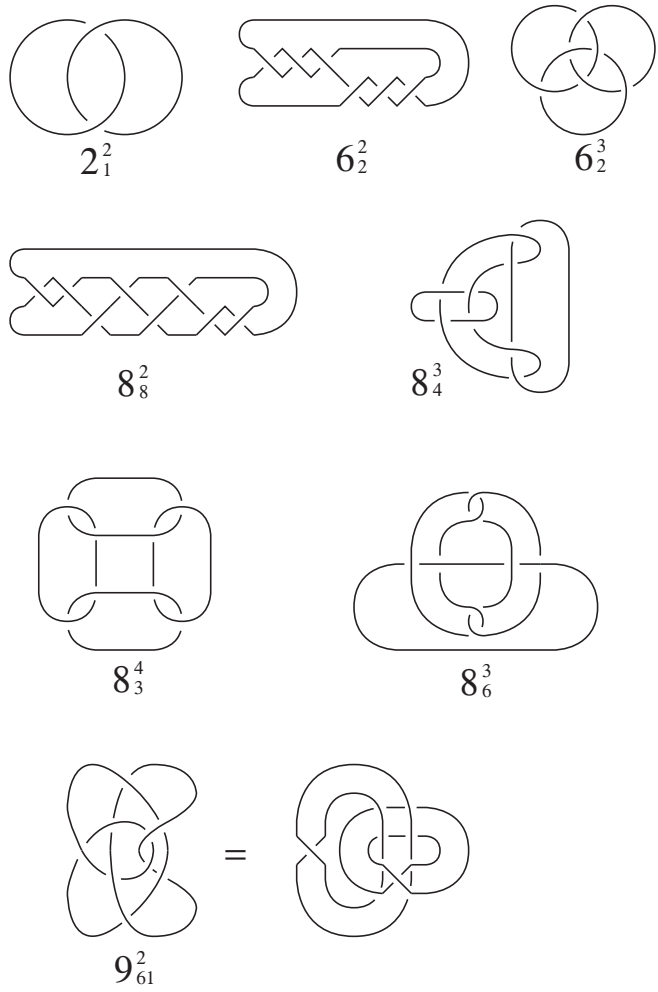} 
\label{fig:9cross}
\caption{Prime amphicheiral links with up to $9$ crossings}
\end{center}
\end{figure}

\newpage

\begin{figure}[htbp]
\begin{center}
\includegraphics[scale=0.5]{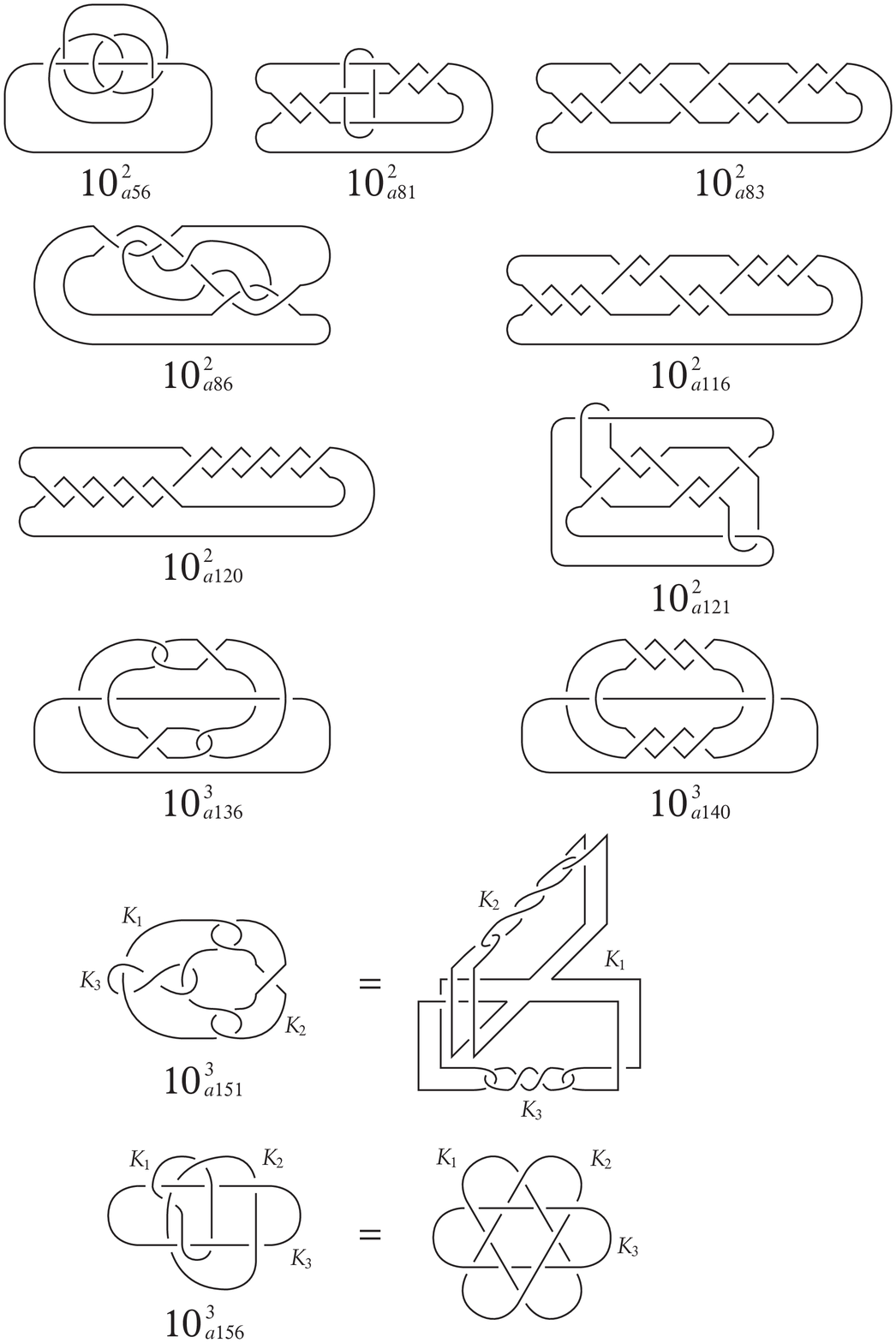} 
\label{fig:11cross1}
\caption{Prime amphicheiral links with $10$ or $11$ crossings 1}
\end{center}
\end{figure}

\newpage

\begin{figure}[htbp]
\begin{center}
\includegraphics[scale=0.5]{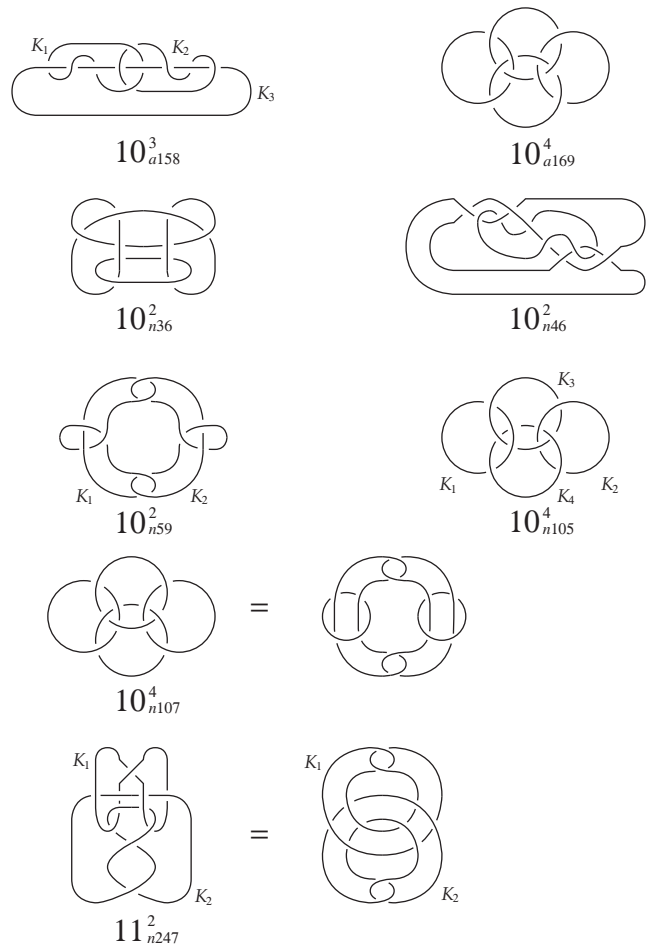} 
\label{fig:11cross2}
\caption{Prime amphicheiral links with $10$ or $11$ crossings 2}
\end{center}
\end{figure}

\newpage

Firstly, we determine which links are amphicheiral,
and whether they are in 
$\mathcal{C}_{11}$, 
$\mathcal{A}_{11}\setminus \mathcal{C}_{11}$, 
$\mathcal{C}_{11}^{\pm}$, or
$\mathcal{A}_{11}^{\pm}\setminus \mathcal{C}_{11}^{\pm}$.
By Figure 3 and Figure 4,
$10_{a56}^2$, $10_{a81}^2$, $10_{a83}^2$, $10_{a86}^2$, 
$10_{a116}^2$, $10_{a120}^2$, 
$10_{a121}^2$, $10_{a136}^3$, $10_{a140}^3$, 
$10_{a151}^3$, $10_{a156}^3$, $10_{a158}^3$, 
$10_{a169}^4$, $10_{n36}^2$, $10_{n46}^2$, 
$10_{n59}^2$, $10_{n105}^4$, $10_{n107}^4$
and $11_{n247}^2$ are amphicheiral.
The reader can check by the figures whether they are in
$\mathcal{C}_{11}$ or $\mathcal{C}_{11}^{\pm}$.
Suppose that $10_{a140}^3$ is $(-, -, -)$-amphicheiral.
Then Milnor's $\bar{\mu}$-invariant $\bar{\mu}(123)$
should be zero.
However since
$${\mit \Delta}_{10_{a140}^3}(t_1, t_2, t_3)
\doteq (t_1-1)(t_2-1)(t_3-1)(t_2t_3+1)^2,$$
we have $\bar{\mu}(123)=\pm 2$ (cf.\ \cite{Co}),
and hence $10_{a140}^3\not \in \mathcal{C}_{11}^-$.
Since $10_{n107}^4$ has a $3$-component sublink
which is equivalent to $6_2^3$,
it cannot be component-preservingly $(-)$-amphicheiral
by Lemma \ref{lm:sub} (1).
However both $10_{a140}^3$ and $10_{n107}^4$
are in $\mathcal{A}_{11}^-\setminus \mathcal{C}_{11}^-$
by suitable orientations.
The $2$-component sublinks of $10_{a151}^3$
are the $2$-component trivial link, 
the positive Whitehead link and the negative Whitehead link.
Since the Whitehead link is not amphicheiral,
$10_{a151}^3\in \mathcal{A}_{11}\setminus \mathcal{C}_{11}$,
and it is in $\mathcal{A}_{11}^+\setminus \mathcal{C}_{11}^+$
and $\mathcal{A}_{11}^-\setminus \mathcal{C}_{11}^-$
by suitable orientations.
Let $L=K_1\cup K_2\cup K_3$ be $10_{a156}^3$ or $10_{a158}^3$.
Then the $2$-component sublinks of $L$
are the $2$-component trivial link $K_1\cup K_2$, 
the positive Hopf link $K_1\cup K_3$ 
and the negative Hopf link $K_2\cup K_3$ by a suitable orientation.
Suppose that $L$ is component-preservingly amphicheiral.
Then 
$${\mit \Delta}_L(t_1, t_2, t_3)
\doteq {\mit \Delta}_L(t_1, t_2, t_3^{-1})$$
should be satisfied by Lemma \ref{lm:amp}.
Since
\begin{eqnarray*}
{\mit \Delta}_{10_{a156}^3}(t_1, t_2, t_3)
& \doteq &
(t_3-1)(t_1t_2t_3-t_1t_2+t_1+t_2-1)
(t_1t_2t_3-t_1t_3-t_2t_3+t_3-1),
\medskip\\
{\mit \Delta}_{10_{a158}^3}(t_1, t_2, t_3)
& \doteq &
(t_3-1)(t_1^2t_2^2+t_1^2t_2t_3-t_1^2t_2-t_1^2t_3+t_1t_2^2t_3
-t_1t_2^2+t_1t_2t_3^2
\medskip\\
& &
-3t_1t_2t_3+t_1t_2-t_1t_3^2+t_1t_3-t_2^2t_3-t_2t_3^2+t_2t_3+t_3^2),
\end{eqnarray*}
$10_{a156}^3$ and $10_{a158}^3$ 
are not component-preservingly amphicheiral.
Suppose that $10_{n105}^4$ is component-preservingly amphicheiral.
Then 
$${\mit \Delta}_{10_{n105}^4}(t_1, t_2, t_3, t_4)
\doteq {\mit \Delta}_{10_{n105}^4}(t_1, t_2, t_3^{-1}, t_4^{-1})$$
should be satisfied by Lemma \ref{lm:amp}.
Since
$${\mit \Delta}_{10_{n105}^4}(t_1, t_2, t_3, t_4)
\doteq
t_1t_2t_3t_4-t_1t_2t_3+t_1t_3^2t_4-t_1t_3t_4-t_2t_3t_4
+t_2t_4-t_3^2t_4+t_3t_4,$$
$10_{n105}^4$ is not component-preservingly amphicheiral.
There is one essential torus $T$ in the exterior of $11_{n247}^2$.
The torus $T$ is trivial as a torus in $S^3$,
and we denote the core of the separated solid torus
by $l_i$\ $(i=1, 2)$
where $l_i$ is in the same connected component of $K_i$.
Suppose that $11_{n247}^2$ is component-preservingly amphicheiral.
Then two links $K_1\cup l_2$ and $K_2\cup l_1$ are amphicheiral.
However since they are the positive Hopf link and
the negative Hopf link, respectively, and they are not amphicheiral,
$11_{n247}^2$ is not component-preservingly amphicheiral.

\medskip

By applying Lemma \ref{lm:Jones} (4) for the rest links, 
we can see non-amphicheirality of them except
$10_{a51}^2$, $10_{a57}^2$, $10_{a171}^4$, 
$10_{n49}^2$, $10_{n93}^3$, $10_{n108}^4$, 
$11_{n127}^2$, $11_{n158}^2$, $11_{n162}^2$, 
$11_{n205}^2$, $11_{n423}^3$, $11_{n432}^3$
and $11_{n437}^3$.
By Lemma \ref{lm:amp}, Lemma \ref{lm:sub},
Lemma \ref{lm:even}, Lemma \ref{lm:2-comp}
and Lemma \ref{lm:sign}
(cf.\ \cite[Example 4.5 and Example 4.6]{KK}),
we can see non-amphicheirality of the links above 
except $10_{a51}^2$ and $11_{n127}^2$.

\begin{figure}[htbp]
\begin{center}
\includegraphics[scale=0.6]{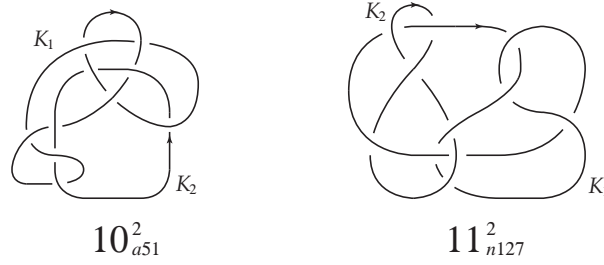} 
\label{fig:1051}
\caption{$10_{a51}^2$ and $11_{n127}^2$}
\end{center}
\end{figure}

\medskip

Let $L=K_1\cup K_2$ be an oriented 
$10_{a51}^2$ or $11_{n127}^2$ as in Figure 5,
and $L'=(-K_1)\cup K_2$.
The HOMFLY polynomials of $10_{a51}^2$ and $11_{n127}^2$ are
\begin{eqnarray*}
P_{10_{a51}^2}(m, l) & = & P_{(10_{a51}^2)'}(m, l^{-1}) \medskip\\
& = & -m^7l^{-1}+m^5(l+4l^{-1}+2l^{-3})
-m^3(2l+6l^{-1}+5l^{-3}+l^{-5})\medskip\\
& & +m(l+2l^{-1}+3l^{-3}+l^{-5})+m^{-1}(l^{-1}+l^{-3}), \medskip\\
P_{11_{n127}^2}(m, l) & = & P_{(11_{n127}^2)'}(m, l^{-1}) \medskip\\
& = & m^5l-m^3(2l^3+4l+l^{-1})
+m(l^5+5l^3+5l+2l^{-1})\medskip\\
& & -m^{-1}(l^5+2l^3+2l+l^{-1}).
\end{eqnarray*}
By Lemma \ref{lm:Jones} (6), we cannot show they are not amphicheiral
by the HOMFLY polynomials.

\medskip

The components of $10_{a51}^2$ consist of
two trivial knots, and the linking number of it is $1$.
Since the Alexander polynomial of it is
$${\mit \Delta}_{10_{a51}^2}(t_1, t_2)\doteq
(t_1t_2^2-2t_1t_2+t_1-t_2^2+t_2-1)
(t_1t_2^2-t_1t_2+t_1-t_2^2+2t_2-1),$$
$10_{a51}^2$ is component-preservingly amphicheiral
by Lemma \ref{lm:amp} if it is amphicheiral.
Since we have
$${\mit \Delta}_{10_{a51}^2}(t, -1)=(4t-3)(3t-4),$$
$10_{a51}^2$ is not amphicheiral by Lemma \ref{lm:branch} (2) and (3).

\medskip

The components of $11_{n127}^2$ consist of
one trivial knot and one figure eight knot, 
and the linking number of it is $1$.
It is component-preservingly amphicheiral
by Lemma \ref{lm:sub} (1) if it is amphicheiral.
Since the Jones polynomial and the HOMFLY polynomial
of $\tilde{K}$ which is the lifted knot of $K_1$
in the $2$-fold branched covering over $K_2$
(cf.\ Figure 6)
are
\begin{eqnarray*}
V_{\tilde{K}}(t) & = &
-2t^6+4t^5-6t^4+9t^3-9t^2+9t-7+5t^{-1}-3t^{-2}+t^{-3},
\medskip\\
P_{\tilde{K}}(m, l) & = &
m^4(2l^2-1)-m^2(4l^4+3l^2-l^{-2})+(2l^6+4l^4+2l^2+1),
\end{eqnarray*}
$11_{n127}^2$ is not amphicheiral
by Lemma \ref{lm:branch} (1), Lemma \ref{lm:Jones} (3) and (6).

\begin{figure}[htbp]
\begin{center}
\includegraphics[scale=0.6]{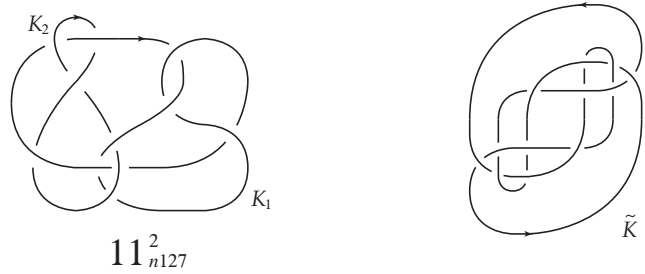} 
\label{fig:1127}
\caption{The lifted knot $\tilde{K}$ of $K_1$ 
in the $2$-fold branched covering over $K_2$}
\end{center}
\end{figure}

We remark that we cannot complete the proof only by
using the Alexander polynomial.
Since the Alexander polynomial of $11_{n127}^2$ is
\begin{equation}\label{eq:11n127}
\begin{matrix}
{\mit \Delta}_{11_{n127}^2}(t_1, t_2) & \doteq &
t_1^2t_2^2-t_1^2t_2+t_1^2
-2t_1t_2^2+t_1t_2-2t_1+t_2^2-t_2+1\hfill\medskip\\
& \doteq & 
\left\{ \left(t_1+t_1^{-1}\right)-2\right\}
\left\{ \left(t_2+t_2^{-1}\right)-1\right\}-1\hfill
\end{matrix}
\end{equation}
we have
$$|{\mit \Delta}_{11_{n127}^2}(-1, -1)|=11\not \equiv 1\quad
(\mathrm{mod}\ \! 4),$$
and hence $11_{n127}^2$ is not $(-, +)$-amphicheiral 
by Corollary \ref{co:branch}.
Let $\zeta_p$ be a primitive $p$-th root of unity for $p\ge 2$.
Then there exists $f(t)\in \mathbb{Z}[t^{\pm 1}]$ such that
$${\mit \Delta}_{\tilde{K}}(t)\doteq
\prod_{i=1}^{p-1}
{\mit \Delta}_{11_{n127}^2}(t, \zeta_p^i)
\doteq
\left\{
\begin{array}{ll}
(f(t))^2(3t^2-5t+3) & (\mbox{$p$ is even}),
\medskip\\
(f(t))^2 & (\mbox{$p$ is odd}),
\end{array}
\right.$$
where
$$f(t)\doteq 
\prod_{1\le i<\frac p2}{\mit \Delta}_{11_{n127}^2}(t, \zeta_p^i)
\in \mathbb{Z}[t^{\pm 1}]\ 
\mbox{(by (\ref{eq:11n127}))}.$$
Hence we cannot prove that $11_{n127}^2$ is not $(+, -)$-amphicheiral 
by Lemma \ref{lm:branch} (2).

\section{Further remarks}\label{sec:remark}
(1)\ In \cite{Kd2}, the author raised a conjecture:

\begin{cj}\label{cj:Kd2}{\rm (\cite[Conjecture 1.1]{Kd2})}
For an $r$-component algebraically split 
component-preservingly amphicheiral link $L$ with $r$ even,
we have
${\mit \Delta}_L(t_1, \ldots, t_r)=0$.
\end{cj}

We gave a partial affirmative answer in \cite[Theorem 1.3]{KK} for the case
that $L$ is an algebraically split component-preservingly 
$(\varepsilon)$-amphicheiral link with even components.
In the prime links with up to $11$ crossings,
only $10_{n36}^2$ and $10_{n107}^4$ are 
algebraically split component-preservingly amphicheiral links
with even components.
They are also component-preservingly $(+)$-amphicheiral links
(i.e.\ $10_{n36}^2, 10_{n107}^4\in \mathcal{C}_{11}^+$).
We can confirm that the Alexander polynomials of them are $0$.
The condition ``component-preservingly" is needed.
$10_{a151}^3$, $10_{n59}^2$ and $11_{n247}^2$
are algebraically split amphicheiral links in 
$\mathcal{A}_{11}\setminus \mathcal{C}_{11}$ and 
$\mathcal{A}_{11}^{\pm}\setminus \mathcal{C}_{11}^{\pm}$
whose Alexander polynomials are
\begin{eqnarray*}
{\mit \Delta}_{10_{a151}^3}(t_1, t_2, t_3) & \doteq &
(t_1-1)(t_2-1)(t_3-1)(t_3^2-3t_3+1)
\medskip\\
{\mit \Delta}_{10_{n59}^2}(t_1, t_2) & \doteq &
(t_1-1)(t_2-1)(t_1-t_2)(t_1t_2-1)
\medskip\\
{\mit \Delta}_{11_{n247}^2}(t_1, t_2) & = &
0.
\end{eqnarray*}
${\mit \Delta}_{10_{n59}^2}(t_1, t_2)$ 
satisfies the condition
$${\mit \Delta}_{10_{n59}^2}(t, t)
={\mit \Delta}_{10_{n59}^2}(t, t^{-1})=0$$
in Lemma \ref{lm:sign}.
We can find examples of $\lambda$-component algebraically split links
in $\mathcal{C}_n\setminus \mathcal{C}_n^{\pm}$
with $\lambda \ge 4$ even in the Milnor links (see Figure 7).
The Alexander polynomials of them are $0$
(cf.\ \cite[Example 6.1 (1)]{Kd2}).

\begin{figure}[htbp]
\begin{center}
\includegraphics[scale=0.55]{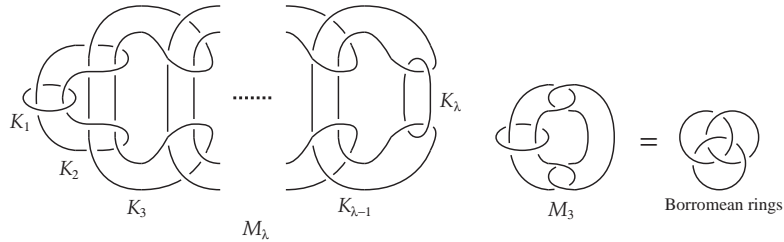} 
\label{fig:Milnor}
\caption{$\lambda$-component Milnor link $M_{\lambda}$}
\end{center}
\end{figure}

\noindent
(2)\ 
In A.~Stoimenow \cite{St1},
the following conjecture is raised as Tait's conjecture IV.

\begin{cj}\label{cj:Tait}{\rm (\cite[Conjecture 2.4]{St1})}
The minimal crossing number of an amphicheiral knot is even.
\end{cj}

If there are counterexamples, then the knots are not alternating
by Lemma \ref{lm:Tait}.
A.~Stoimenow \cite{St2} found prime amphicheiral knots
with the odd minimal crossing number $c$ for every $c\ge 15$.
He points out that if $c\le 13$, then Conjecture \ref{cj:Tait}
is affirmative.
We have already found counterexamples
$9_{61}^2$ and $11_{n247}^2$ for the case of links.
However they are not component-preservingly amphicheiral.
Recently Y.~Kobatake found an example of 
a $2$-component prime component-preservingly amphicheiral link 
with the minimal crossing number $21$ and with the linking number $3$,
whose components consist of the unknot and
an amphicheiral knot $15_{224980}$ in the table of \cite{HTW}.

\bigskip

{\noindent {\bf Acknowledgements}}\ 
The author would like to express gratitude to
Akio Kawauchi for giving him useful comments.

{\footnotesize
\par
\medskip
Teruhisa KADOKAMI\par 
Department of Mathematics,\par
East China Normal University,\par
Dongchuan-lu 500, Shanghai, 200241, China \par
{\tt mshj@math.ecnu.edu.cn}\par
{\tt kadokami2007@yahoo.co.jp}\par
}
\end{document}